\documentclass{tac}

\usepackage{enumerate}
\usepackage{amsmath,xspace,amssymb,mathrsfs}

\input xy
\xyoption{all}
\CompileMatrices

\title{An embedding theorem for adhesive categories}
\author{Stephen Lack} 
\thanks{The support of the Australian Research Council and
DETYA is gratefully acknowledged.} 
\address{Mathematics Department NSW 2109 \\
Australia}
\eaddress{steve.lack@mq.edu.au}
\copyrightyear{2011}

\newcommand{\B}{{\ensuremath{\mathscr B}}\xspace}
\newcommand{\C}{{\ensuremath{\mathscr C}}\xspace}
\newcommand{\E}{{\ensuremath{\mathscr E}}\xspace}
\newcommand{\op}{\ensuremath{{}^{\textrm{op}}}\xspace}
\newcommand{\Set}{\textnormal{\bf Set}\xspace}

\DeclareMathOperator\Span{Span}
\DeclareMathOperator\Sh{Sh}

\def\x{\times}

\begin{document}

\label{firstpage}
\maketitle

\begin{abstract}
Adhesive categories are categories which have pushouts with one leg a monomorphism, all pullbacks, and certain exactness conditions relating these pushouts and pullbacks. We give a new proof of the fact that every topos is adhesive. We also prove a converse: every small adhesive category has a fully faithful functor in a topos, with the functor preserving the all the structure. Combining these two results, we see that the exactness conditions in the definition of adhesive category are exactly the relationship between pushouts along monomorphisms and pullbacks which hold in any topos.
\end{abstract}

\section{Introduction}

Many different categorical structures involve certain limits and colimits
connected by exactness conditions which state, roughly speaking, that 
limits and colimits in the category interact in the same way that they do in 
the category of sets; or better perhaps, that they interact in the same way that
they do in any topos. 

For example, we have the structure of regular category, which can be characterized
as the existence of finite limits and coequalizers of kernel pairs, along with the 
condition that these coequalizers are stable under pullback. As evidence that
this condition characterizes the interaction between finite limits and coequalizers
of kernel pairs in a topos, we have on the one hand, the fact that any topos is a 
regular category, and on the other hand, the theorem of Barr \cite{Barr:exact}
which asserts that  any small regular category has a fully faithful embedding
in a topos, and this embedding preserves coequalizers of kernel pairs and finite limits.

Then again, there is the structure of extensive category, which can be characterized
as the existence of finite coproducts and pullbacks along coproduct injections,
along with the condition that these coproducts are stable under pullback and 
disjoint. Once again, every topos is extensive, and every small extensive 
category has a fully faithful embedding in a topos, and this embedding preserves 
finite coproducts and all existing (finite) limits \cite{Lawvere-Como}.

The notion of adhesive category was introduced in \cite{adh} as a categorical
framework for graph transformation and rewriting; see the volume 
\cite{GraphTransformations} for more on this point of view. Adhesive categories 
are the analogue of extensive categories where one works with pushouts along a 
monomorphism rather than coproducts; the definition is recalled in the following
section. Once again, every topos is adhesive, as
shown in \cite{topos-adh} and in Theorem~\ref{thm:topos-adh} below, 
and it is natural to ask 
whether there is a corresponding embedding theorem.
The purpose of this paper is to show that this is the 
case: every small adhesive category has a fully faithful embedding in a topos, 
and this embedding preserves pushouts along monomorphisms and all
existing finite limits. I am grateful to Bill Lawvere for suggesting the question.


In each case the construction is essentially the same. The Yoneda embedding
$Y:\C\to[\C\op,\Set]$ preserves all existing limits, but few colimits; in order
to correct this we replace $[\C\op,\Set]$ by the full subcategory of sheaves
for a topology on \C chosen so that the representables are sheaves and the
restricted Yoneda embedding $\C\to\Sh\C$ preserves the colimits in question.
A general analysis of the sorts of colimits and exactness properties
that can be dealt with in this way will be given in \cite{GarnerLack}.

From a model-theoretic point of view, these embedding theorems can be viewed as 
completeness results. In the case of adhesive categories, for instance, since
every topos is adhesive, and every (small) adhesive category has a fully faithful
embedding into a topos preserving pushouts along monomorphisms and pullbacks,
the adhesive category axioms capture that fragment of the structure of a topos
which involves pushouts along monomorphisms and pullbacks. 
See \cite{MakkaiPare} for a detailed treatment of this point of view.


\section{Adhesive categories}

Recall that a category \C with finite coproducts is {\em extensive} if for
all objects $A$ and $B$ the coproduct functor $\C/A\x\C/B\to\C/(A+B)$ 
is an equivalence of categories. This is equivalent to saying
that \C has finite coproducts and pullbacks along coproduct injections, and 
in a commutative diagram
$$\xymatrix{
A' \ar[r] \ar[d] & E \ar[d] & B' \ar[l] \ar[d] \\
A \ar[r] & A+B & B \ar[l] }$$
in which the bottom row is a coproduct,  the top row is a coproduct if and
only if the squares are pullbacks.

In a pushout square 
$$\xymatrix{ 
C \ar[r]^g \ar[d]_f & B \ar[d] \\
A \ar[r] & D }$$
the role of the morphisms $f$ and $g$ is entirely symmetric, but sometimes
we wish to think of it non-symmetrically, and call it a {\em pushout along $f$},
or alternatively a {\em pushout along $g$}. If either $f$ or $g$ is a monomorphism,
we call it a {\em pushout along a monomorphism}. 

A category with pullbacks is said to be {\em adhesive} if it has pushouts
along monomorphisms, every such pushout is stable under pullback, and
in a cube
$$\xymatrix @R1pc @C1pc {
& C' \ar[rr] \ar[dl] \ar[dd] && B' \ar[dl] \ar[dd] \\
A' \ar[rr] \ar[dd] && D' \ar[dd] \\
& C \ar[rr] \ar[dl] && B \ar[dl] \\
A \ar[rr] && D }$$
in which the top and bottom squares are pushouts along monomorphisms, and the 
back and left squares are pullbacks, then the remaining squares are
pullbacks. Alternatively, one can combine the stability condition and the 
other condition into the single requirement that in such a cube, if the bottom
square is a pushout along a monomorphism, and the left and back squares are
pullbacks, then the top square is a pushout if and only if the front and right 
squares are pullbacks.

Various simple facts follow; for example the following are both proved in \cite{adh}:


\begin{proposition}\label{prop:basic}
  In any adhesive category, the pushout of a monomorphism along any map 
is a monomorphism, and the resulting square is also a pullback.
\end{proposition}

\proof
Suppose that 
$$\xymatrix @R1pc @C1pc {
& C \ar[rr] \ar[dl]_m && B \ar[dl]^{n} \\
A \ar[rr] && D }$$
is a pushout in which $m$ is a monomorphism. We have to show that $n$ is
a monomorphism and that the square is also a pullback. 
Expand the given square into a cube, as in the diagram 
$$\xymatrix @R1pc @C1pc {
& C \ar[rr] \ar@{=}[dl] \ar@{=}[dd] && B \ar@{=}[dl] \ar@{=}[dd] \\
C \ar[rr] \ar[dd] && B \ar[dd] \\
& C \ar[rr] \ar[dl] && B \ar[dl] \\
A \ar[rr] && D }$$
in which the front and bottom faces are both the original square.
The top face is a pushout, the back face a pullback, the left
also a pullback (because $m$ is a monomorphism), and so the front and right
faces are pullbacks. The front face being a pullback means that the original
square is a pullback; the right face being a pullback means that $n$ is 
a monomorphism.
\endproof

The proof of the following result is slightly more complicated; it can 
be found in \cite{qadh}. 

\begin{proposition}\label{prop:union}
An adhesive category has binary unions of subobjects, and they are effective.
\end{proposition}

The claim that the unions are effective means that the union of a pair of subobjects
is constructed as the pushout over their intersection.

There is a variant of adhesive categories called quadiadhesive categories,
which uses pushouts along regular monomorphisms rather than along all
mono\-morphisms. Perhaps suprisingly, and contrary to what is claimed in 
\cite{qadh}, it is not the case that every quasitopos is quasiadhesive: see
\cite{JohnstoneLackSobocinski}.

An elegant reformulation of the adhesive condition was given in 
\cite{HeindelSobocinski} using the 
bicategory $\Span(\C)$ of spans in \C. This bicategory has the same objects as \C,
while a morphism from $A$ to $B$ is a diagram $A \gets E \to B$, with composition
given by pullback. There is an inclusion pseudofunctor $\C\to\Span(\C)$
sending a morphism $f:A\to B$ to the span from $A=A\to B$ with left leg
the identity on $A$ and right leg just $f$. The reformulation then says that
a category \C with pullbacks is adhesive if and only if it has pushouts along
monomorphisms, and these are sent by the inclusion $\C\to\Span(\C)$ to {\em bicolimits} in $\Span(\C)$.

A functor between adhesive categories is called adhesive if it preserves
pushouts along monomorphisms and pullbacks.

\section{Adhesive categories and toposes}

In this section we see that every topos is adhesive, and so every full subcategory
of a topos, closed under pushouts along monomorphisms and pullbacks, is 
adhesive. Then we see the converse: every (small) adhesive category has a fully
faithful adhesive functor into a topos. 

First we give a proof, different to that of \cite{topos-adh}, that every topos
is adhesive. It uses one of the Freyd embedding theorems to reduce
to the case of a Boolean topos. It is also possible to give direct proofs, 
but these all seem to be rather longer --- see \cite{topos-adh} for 
one possibility.

\begin{theorem}\label{thm:topos-adh}
Every topos \E is adhesive. 
\end{theorem}

\proof
Recall that a topos  is  Boolean when  every subobject is 
complemented, in the sense that it is a coproduct injection. By a 
theorem of Freyd~\cite{Freyd-aspects},  for every topos \E there is a
Boolean topos \B and a faithful  functor from \E to \B
which preserves finite limits and finite (in fact any) colimits. Since the 
functor is faithful and preserves finite limits and colimits, it also reflects
isomorphisms, and so it  reflects finite limits and colimits. 
Thus the adhesive category axioms for \E will follow from those for \B,
and it suffices to prove that any Boolean topos is adhesive.

In a Boolean topos, any monomorphism $m:C\to A$ is a coproduct 
injection $m:C\to C+X$, and now for any map $f:C\to B$ the corresponding
pushout is the square
$$\xymatrix{
C \ar[d]_{m} \ar[r]^{f} & B \ar[d]^{n} \\
C+X \ar[r]_{f+X} & B+X }$$
with $n:B\to B+X$ once again the coproduct injection.

Now consider a cube over this square, with back and left faces pullbacks. 
By extensivity, this has the form 
$$\xymatrix @R1pc @C1pc {
& C' \ar[rr]^{f'} \ar[dl]_{m'} \ar[dd]^(0.3)c && B' \ar[dl]_{n'} \ar[dd]^b \\
C'+X' \ar[rr]^(0.7){g'} \ar[dd]_{c+x} && D \ar[dd]^(0.2){d} \\
& C \ar[rr]^(0.7)f \ar[dl]_m && B \ar[dl]^n \\
C+X \ar[rr]_{f+X} && B+X }$$
in which $m':C'\to C'+X'$ is the coproduct injection.

Of course if the front and right faces are pullbacks, then the top face is 
a pushout, by stability of colimits under pullback. Conversely, if the top face is a 
pushout, then it has the form 
$$\xymatrix{
C' \ar[r]^{f'} \ar[d]_{m'} & B' \ar[d]^{n'} \\
C'+X' \ar[r]_{f'+X'} & B'+X' }$$
and now $d$ is just $b+x:B'+X'\to B+X$ by commutativity of the 
front and right faces. The fact that the front and right
faces are pullbacks then follows by extensivity.
\endproof

Before we prove the embedding theorem, we need the following lemma.


\begin{lemma}\label{lemma:basic}
Let a monomorphism $m:C\to A$ and a map $f:C\to B$ be given in any  
adhesive category \C, and construct the diagram
$$\xymatrix{
C \ar[r]^{\gamma} \ar[d]_m & C_2 \ar[d]_{m_2} \ar@<1ex>[r]^{f_1}\ar@<-1ex>[r]_{f_2}
& C \ar[r]^{f} \ar[d]_m & B \ar[d]^n \\
A \ar[r]_{\delta} & A_2 \ar@<1ex>[r]^{g_1} \ar@<-1ex>[r]_{g_2} & A \ar[r]_g & D }$$
in  which the right hand square is a pushout, $C_2$ and $A_2$ are the
kernel pairs of $f$ and $g$, and $\gamma$ and $\delta$ are the diagonal maps.
Then the left hand square is a pushout and a pullback.
\end{lemma}

\proof
Pulling back the right hand square along $g:A\to D$ gives the square
$$\xymatrix{
C_2 \ar[r]^{f_1} \ar[d]_{m_2} & C \ar[d]^{m} \\
A_2 \ar[r]_{g_1} & A }$$
which is therefore a pushout. Form the diagram
$$\xymatrix{
C \ar[r]^\gamma \ar[d]_m & C_2 \ar[r]^{f_1} \ar[d]^j & C \ar[d]^m \\
A \ar[r]^i \ar@{.>}[dr]_\delta & E\ar[r]^{h_1} \ar[d]^k & A \ar[d]^1 \\
& A_2 \ar[r]_{g_1} & A }$$
in which the top left square is a pushout, $h_1$ is the unique morphism satisfying
$h_1i=1$ and $h_1j=mf_1$, and 
$k$ is the unique morphism satisfying $ki=\delta$ and $kj=m_2$. By
Proposition~\ref{prop:union}, the morphism $k:E\to A_2$ is the union of
the subobjects $\delta:A\to A_2$ and $m_2:C_2\to A_2$, and so in particular
is a monomorphism. The lemma asserts that it is invertible. 

Now the top left square and the composite of the upper squares are 
both pushouts, so the top right square is also a pushout, by the cancellativity
properties of pushouts. The composite of the two squares on the right is the pushout constructed at the beginning of the proof,
so finally the lower right square is a pushout by the cancellativity property of 
pushouts once again.

Since $k$ is a monomorphism, the lower right square is a pullback by 
 Proposition~\ref{prop:basic}. Thus $k$ is invertible, and our square is indeed
a pushout. It is a pullback by Proposition~\ref{prop:basic} again along with the fact
that $m$ is a monomorphism. 
\endproof


\begin{theorem}
Any small adhesive category admits a full adhesive embedding into a topos.
\end{theorem}

\proof
Let \C be a small adhesive category. Consider the topology generated
by all pairs $(g,n)$ arising as above in a pushout
$$\xymatrix{
C \ar[r]^{f} \ar[d]_{m} & B \ar[d]^{n} \\
A \ar[r]_{g} & D}$$
in which $m$ is a monomorphism. The pushout square is also a pullback
by Proposition~\ref{prop:basic}, thus $f$ and $m$ are determined up to 
isomorphism by $g$ and $n$.
These $(g,n)$  generate a topology since pushouts
along monomorphisms are stable under pullback, thus the category \E of
sheaves for the topology is a topos. We shall show that the Yoneda embedding
lands in \E, giving a fully faithful functor $Y:\C\to\E$, and that $Y$ is adhesive.

To do this, we  need to understand the sheaf condition. 
Consider a pair $(g,n)$ arising as above.  Since $n$ is a monomorphism, by 
Proposition~\ref{prop:basic} once again, the kernel pair of $n$ is just $B$.
Let $g_1,g_2:A_2\rightrightarrows A$ be the kernel pair of $g$. Then a functor 
$F:\C\op\to\Set$ satisfies the sheaf condition precisely when, for 
each such $(g,n)$, the morphisms $Fg$ and $Fn$ exhibit $FD$ as the
limit of the diagram 
$$\xymatrix{
& FC & FB \ar[l]_{Ff} \\
FA_2 & FA \ar@<1ex>[l]^{Fg_1} \ar@<-1ex>[l]_{Fg_2} \ar[u]^{Fm} &
FD \ar[l]_{Fg} \ar[u]_{Fn} }$$
in \Set. 

Clearly the representables satisfy this condition, so the Yoneda
embedding lands in \E, giving a fully faithful functor $Y:\C\to\E$,
which preserves all existing limits, and in particular all finite limits.
We must show that it also preserves pushouts along monomoprhisms.
This is equivalent to the condition that every sheaf $F$ sends pushouts
along monomorphisms in $\C$ to limits in \Set. If $g$ were monomorphic,
so that the kernel pair $A_2$ of $g$ were just $A$, this would be immediate,
but in general there is a little work to do.

We must show that if $x\in FA$ and $y\in FB$ satisfy $Fm.x=Ff.y$, then also 
$Fg_1.x=Fg_2.x$, so that by the sheaf condition $x$ and $y$ arise from some
(necessarily unique) $z\in FD$. 

To do this, we use Lemma~\ref{lemma:basic}. This tells us that 
$A_2$ is the pushout of $A$ and $C_2$ over $C$, and 
all maps in this pushout are monomorphisms, so by the sheaf condition
$$\xymatrix{
FC & FC_2 \ar[l]_{F\gamma} \\
FA \ar[u]^{Fm} & FA_2 \ar[l]^{F\delta} \ar[u]_{Fm_2} }$$
is a pullback, and so $F\delta$ and $Fm_2$ are jointly monic. Now
$$F\delta.Fg_1.x=F(g_1\delta).x=x=F(g_1\delta).x=F\delta.Fg_2.x$$
$$Fm_2.Fg_1.x=Ff_1.Fm.x=Ff_1.Ff.y=Ff_2.Ff.y=Ff_2.Fm.x=Fm_2.Fg_2.x$$
and so $Fg_1.x=Fg_2.x$
as required.
\endproof



\begin{thebibliography}{1}

\bibitem{GraphTransformations}
Andrea~Corradini et~al, editor.
\newblock {\em Graph Transformations}, volume 4178 of {\em Lecture Notes in
  Comput. Sci.}
\newblock Springer, 2006.

\bibitem{Barr:exact}
Michael Barr.
\newblock Exact categories and categories of sheaves.
\newblock {\em Lecture Notes Math.} 236:1--120, Springer, Berlin, 1971.



\bibitem{Freyd-aspects}
Peter Freyd.
\newblock Aspect of topoi.
\newblock {\em Bull. Austral. Math. Soc.}, 7:1--76, 1972.

\bibitem{GarnerLack}
Richard Garner and Stephen Lack.
\newblock Lex colimits, in preparation.

\bibitem{HeindelSobocinski}
Tobias Heindel and Pawe{\l} Soboci{\'n}ski.
\newblock Van {K}ampen colimits as bicolimits in span.
\newblock In {\em Algebra and coalgebra in computer science}, volume 5728 of
  {\em Lecture Notes in Comput. Sci.}, pages 335--349. Springer, Berlin, 2009.

\bibitem{JohnstoneLackSobocinski}
Peter~T. Johnstone, Stephen Lack, and Pawe{\l} Soboci{\'n}ski.
\newblock Quasitoposes, quasiadhesive categories and {A}rtin glueing.
\newblock In {\em Algebra and coalgebra in computer science}, volume 4624 of
  {\em Lecture Notes in Comput. Sci.}, pages 312--326. Springer, Berlin, 2007.

\bibitem{adh}
Stephen Lack and Pawe{\l} Soboci{\'n}ski.
\newblock Adhesive categories.
\newblock In {\em Foundations of software science and computation structures},
  volume 2987 of {\em Lecture Notes in Comput. Sci.}, pages 273--288. Springer,
  Berlin, 2004.

\bibitem{qadh}
Stephen Lack and Pawe{\l} Soboci{\'n}ski.
\newblock Adhesive and quasiadhesive categories.
\newblock {\em Theor. Inform. Appl.}, 39(3):511--545, 2005.

\bibitem{topos-adh}
Stephen Lack and Pawe{\l} Soboci{\'n}ski.
\newblock Toposes are adhesive.
\newblock In {\em Graph transformations}, volume 4178 of {\em Lecture Notes in
  Comput. Sci.}, pages 184--198. Springer, Berlin, 2006.

\bibitem{Lawvere-Como}
F.~William Lawvere.
\newblock Some thoughts on the future of category theory.
\newblock In {\em Category theory ({C}omo, 1990)}, volume 1488 of {\em Lecture
  Notes in Math.}, pages 1--13. Springer, Berlin, 1991.

\bibitem{MakkaiPare}
Michael Makkai and Robert Par{\'e}.
\newblock {\em Accessible categories: the foundations of categorical model
  theory}, volume 104 of {\em Contemporary Mathematics}.
\newblock American Mathematical Society, Providence, RI, 1989.

\end{thebibliography}


\end{document}